\documentclass{article}
\pdfoutput=1

\usepackage{fontenc}
\usepackage{graphicx, float}
\usepackage[a4paper, total={6in, 9.5in}]{geometry}
\usepackage{parskip}

\usepackage[fleqn,leqno]{amsmath}
\usepackage{amsthm, amsfonts, amssymb}
\usepackage{thmtools}

\usepackage{hyperref, bookmark}
\hypersetup{hypertexnames = false, bookmarksdepth = 2, bookmarksopen = true, colorlinks, linkcolor = black, citecolor = black, urlcolor = black, pdfstartview={XYZ null null 1}}
\usepackage[capitalise]{cleveref}
\usepackage{colonequals}

\usepackage{libertine}
\usepackage[libertine]{newtxmath}
\usepackage{eucal}
\usepackage{microtype}

\usepackage{booktabs}

\usepackage{tikz, tikz-cd}
\usetikzlibrary{calc}
\usetikzlibrary{arrows, matrix}
\usetikzlibrary{positioning}

\usepackage[backend=biber, datamodel=mrnumber, maxbibnames=99, sortcites]{biblatex}
\declaretheoremstyle[
  spaceabove = 3pt,
  spacebelow = 3pt,
  bodyfont = \itshape,
]{first}
\declaretheoremstyle[
  spaceabove = 3pt,
  spacebelow = 3pt,
]{second}
\declaretheorem[numberwithin=section, style=first]{theorem}
\declaretheorem[sibling=theorem, style=first]{conjecture}
\declaretheorem[sibling=theorem, style=first]{corollary}
\declaretheorem[sibling=theorem, style=first]{lemma}
\declaretheorem[sibling=theorem, style=first]{proposition}

\declaretheorem[sibling=theorem, style=second]{example}
\declaretheorem[sibling=theorem, style=second]{remark}
\declaretheorem[sibling=theorem, style=second]{definition}

\declaretheorem[sibling=theorem, style=second]{assumption}

\Crefname{assumption}{Assumption}{Assumptions}
\Crefname{convention}{Convention}{Conventions}
\Crefname{setup}{Setup}{Setups}

\declaretheorem[numberwithin=section, style=first, title=Theorem]{alphatheorem}

\declaretheorem[sibling=alphatheorem, style=first, title=Question]{alphaquestion}

\crefname{alphatheorem}{Theorem}{Theorems}
\crefname{alphaconjecture}{Conjecture}{Conjectures}
\crefname{alphaquestion}{Question}{Questions}

\DeclareMathOperator{\Hom}{Hom}
\DeclareMathOperator{\Ext}{Ext}
\DeclareMathOperator{\ext}{ext}
\DeclareMathOperator{\GL}{GL}
\DeclareMathOperator{\parameterspace}{R}
\DeclareMathOperator{\modulispace}{M}
\DeclareMathOperator{\Pic}{Pic}
\DeclareMathOperator{\rk}{rk}
\DeclareMathOperator{\Hochschild}{HH}
\DeclareMathOperator{\HH}{H}
\DeclareMathOperator{\hh}{h}

\newcommand{\NN}{\mathbb{N}}
\newcommand{\ZZ}{\mathbb{Z}}
\newcommand{\QQ}{\mathbb{Q}}

\newcommand{\PP}{\mathbb{P}}

\DeclareMathOperator{\sheafHom}{\mathscr{H}\kern -2.5pt\mathit{om}}

\mathchardef\mhyphen="2D
\newcommand{\stable}{\mhyphen\mathrm{st}}
\newcommand{\semistable}{\mhyphen\mathrm{sst}}

\newcommand{\tuple}[1]{\ensuremath{\mathbf{#1}}}
\newcommand{\field}{{k}}

\newcommand{\derived}{\mathbf{D}}
\newcommand{\bounded}{\mathrm{b}}
\newcommand{\geomquotient}{\ensuremath{/\kern -2pt/}}
\newcommand{\canonicalstability}{\ensuremath{{\theta_{\mathrm{can}}}}}

\title{Partial semiorthogonal decompositions for quiver moduli}
\author{Gianni Petrella}
\date{April 21, 2025}

\begin{document}

\maketitle
\begin{abstract}
We embed several copies of the derived category of a quiver and certain line bundles
in the derived category of an associated moduli space of representations,
giving the start of a semiorthogonal decomposition.
This mirrors the semiorthogonal decompositions
of moduli of vector bundles on curves.
Our results are obtained with \textsc{QuiverTools},
an open-source package of tools for quiver representations,
their moduli spaces and their geometrical properties.
\end{abstract}
\tableofcontents
\section{Introduction}\label{section:introduction}
There exists a rich dictionary between moduli spaces
of vector bundles on curves and moduli spaces of representations of quivers,
see e.g. \cite{MR3882963}.
Following this dictionary, it is often fruitful to consider an already established
result in the setting of curves and try to prove its appropriate translation
in the case of quivers, or vice versa.
For moduli spaces of vector bundles on curves,
semiorthogonal decompositions of derived categories
are a topic of active interest \cite{MR3713871, MR3880395, MR3764066, MR3954042, MR4651618}.
Following the quiver-curve dictionary,
we can expect the existence of analogous semiorthogonal decompositions
for derived categories of moduli of representations of quivers, for which a first step
has been obtained in \cite{vector-fields-paper}.

We propose several precise questions in this direction in
\cref{question:first-semiorthogonal-embedding,question:r-semiorthogonal-embeddings},
conjecturing the existence of a partial semiorthogonal decomposition
which mirrors the state of the art in the case of curves.
We will subsequently give partial answers to these questions in
several meaningful cases, and a sufficient criterion settling a weaker version
of~\cref{question:r-semiorthogonal-embeddings,question:strongly-exceptional}.
Our proofs make an essential use of three ingredients: explicit computations
in the Chow ring of quiver moduli, Teleman quantization and Serre duality.

The computations in this work are performed using newly added functionalities of
\textsc{QuiverTools} \cite{quivertools}.
It is a software package which deals
with quiver representations and their moduli spaces.
It enables us to reduce complex geometric problems
involving the bounded derived category of coherent sheaves on quiver moduli
to concrete very explicit computations, and to actually perform them.

Let us elaborate about some instances of the aforementioned quiver-curve dictionary
and on how they motivate the questions raised in this paper.

\paragraph{The case of curves}\label{paragraph:SOD-curves}\phantom{a}

Given a curve $C$ of genus $g \geq 2$, a rank $r \geq 2$ and a line bundle $\mathcal{L}$
with $\gcd(r,\deg(\mathcal{L})) = 1$,
there exists the fine moduli space $\modulispace_C(r,\mathcal{L})$ of stable vector bundles
of rank $r$ on $C$ whose determinant is $\mathcal{L}$.
This is a smooth projective Fano variety of
dimension $(r^2 - 1)(g - 1)$, index $2$ and Picard rank $1$.
On the product $C \times \modulispace_C(r,\mathcal{L})$ there exists the so-called
\emph{Poincaré vector bundle} $\mathcal{W}$, which is unique up to a choice of linearisation.

First \cite{MR3713871, MR3880395} and \cite{MR3764066} for the case
$r = 2, \deg(\mathcal{L}) = 1$, followed by \cite{MR3954042}
for $r \geq 2, \deg(\mathcal{L}) = 1$ and by \cite{MR4651618}
for the general case of $\gcd(r, \deg(\mathcal{L})) = 1$,
showed that the Fourier--Mukai transform $\Phi_{\mathcal{W}}$
associated to $\mathcal{W}$ is a fully faithful functor,
and thus embeds $\derived^\bounded(C)$ in $\derived^\bounded(\modulispace_C(r,\mathcal{L}))$.
It is also possible to twist the image of $\Phi_{\mathcal{W}}$ by
$\Theta$, the ample generator of the Picard group,
and embed two copies of $\derived^\bounded(C)$.

In~\cite[Theorem 1.4]{MR4651618}
these copies are shown to be semiorthogonal for $g \geq 6$,
giving rise to a partial semiorthogonal decomposition of $\modulispace_C(r,\mathcal{L})$
of the form
\begin{equation}\label{equation:partial-SOD-VBAC}
\derived^\bounded(\modulispace_C(r,\mathcal{L})) =
\langle \mathcal{O}_{\modulispace}, \Phi_{\mathcal{W}}(\derived^\bounded(C)),
\Theta, \Phi_{\mathcal{W}}(\derived^\bounded(C)) \otimes \Theta;
\ \dots  \rangle.
\end{equation}

The description of a finer semiorthogonal decomposition
of $\derived^\bounded(\modulispace_C(2,\mathcal{L}))$
was conjectured independently in \cite{BGM-conjecture, MR4557892}
and \cite{lee2018remarks}, as the so-called BGMN conjecture.
This not only uses embedded copies of $\derived^\bounded(C)$, but also
$\derived^{\bounded}(\mathrm{Sym}^i(C))$, for $i = 1, \dots, g-1$.
The BGMN conjecture was proven to be true in \cite{MR4850959}.
What the appropriate translation of
this result should be in the context of quiver moduli is unclear, but in this work
we attempt to make the first steps in this direction.

For vector bundles of rank $3$ and degree $1$, Gómez--Lee \cite[Conjecture 1.9]{gómez2020motivicdecompositionsmodulispaces}
propose a similar decomposition involving products of symmetric powers of $C$, while
for $r \geq 4$ there is no precise conjecture - let alone result,
about a finer semiorthogonal decomposition of $\derived^{\bounded}(\modulispace_C(r,\mathcal{L}))$.

\paragraph{The case of quivers}\label{paragraph:SOD-quiver-moduli}\phantom{a}

Let now $Q$, $\tuple{d}$ and $\theta$ be
an acyclic quiver, a dimension vector and a stability parameter respectively.
The moduli space~$\modulispace^{\theta\stable}(Q,\tuple{d})$
of $\theta$-stable representations of $Q$ of dimension vector $\tuple{d}$
exists and is described e.g. in \cite{MR2484736}.
We denote it by $\modulispace$.

Following the dictionary between vector bundles on curves and quiver representations,
\cite{vector-fields-paper} proves that under \cref{assumptions},
there exists a copy of $\derived^\bounded(Q)$ embedded
in $\derived^\bounded(\modulispace)$ via a fully faithful
Fourier--Mukai functor~$\Phi_{\mathcal{U}}$.

Let now $r$ be the index of $\modulispace$, that is, the largest integer that divides
the canonical divisor $K_{\modulispace}$.
Let $H \colonequals - \frac{1}{r}K_{\modulispace}$.
Following the aforementioned dictionary,
we investigate the translation of \eqref{equation:partial-SOD-VBAC}
and compose $\Phi_{\mathcal{U}}$ with the twist by various powers of $\mathcal{O}(H)$,
to embed multiple copies of $\derived^\bounded(Q)$ in $\derived^\bounded(\modulispace)$.
We denote the composition of $\Phi_{\mathcal{U}}$ with the twist
by $\mathcal{O}(H)$ by $\Phi_{\mathcal{U}(H)}$, as in \eqref{equation:twisted-FMT},
and we investigate whether the images of these twisted functors are
semiorthogonal among them and with respect to powers of $\mathcal{O}(H)$.
This leads us to formulate the following questions.

\begin{alphaquestion}\label{question:first-semiorthogonal-embedding}
Let $Q, \tuple{d}$ and $\theta$ be as in \cref{assumptions}.
Is there a semiorthogonal decomposition of the form
\begin{equation}
    \label{equation:first-embedding}
    \derived^{\bounded}(\modulispace) =
    \left\langle \mathcal{O}_{\modulispace}, \Phi_{\mathcal{U}}(\derived^{\bounded}(Q));\ \dots \right\rangle
\end{equation}
for $\derived^{\bounded}(\modulispace) \colonequals \derived^{\bounded}(\modulispace^{\theta\stable}(Q,\tuple{d}))$?
\end{alphaquestion}

\begin{alphaquestion}\label{question:r-semiorthogonal-embeddings}
Assuming a positive answer to \cref{question:first-semiorthogonal-embedding}, if $r \geq 2$,
is there a finer semiorthogonal decomposition of the form
\begin{equation}
    \label{equation:partial-SOD}
    \derived^{\bounded}(\modulispace) =
    \left\langle \mathcal{O}_{\modulispace}, \Phi_{\mathcal{U}}(\derived^{\bounded}(Q)),
    \mathcal{O}(H), \Phi_{\mathcal{U}(H)}(\derived^{\bounded}(Q)),\dots,
    \mathcal{O}((r-1)H), \Phi_{\mathcal{U}((r-1)H)}(\derived^{\bounded}(Q)); \dots \right\rangle
\end{equation}
for $\derived^{\bounded}(\modulispace) \colonequals \derived^{\bounded}(\modulispace^{\theta\stable}(Q,\tuple{d}))$?
\end{alphaquestion}

\begin{remark}\label{remark:existence-of-linearisation}
The universal family on $\modulispace$, and therefore the functor
$\Phi_{\mathcal{U}}$, is again only unique up to the choice of a linearisation,
see \cref{remark:choice-linearisation}.
The answers to
\cref{question:first-semiorthogonal-embedding,question:r-semiorthogonal-embeddings}
depend thus on the chosen linearisation.
Accordingly, the questions
ask for the existence of a choice
for which \cref{equation:first-embedding} and \cref{equation:partial-SOD} hold.
\end{remark}

\begin{remark}
The decomposition in \eqref{equation:partial-SOD} has the shape of a
\emph{Lefschetz decomposition}, see \cite[Definition 2.1]{MR3728631}.
We will not discuss this further in this work.
\end{remark}

\begin{remark}\label{remark:exceptional-collection}
The derived category $\derived^\bounded(Q)$ itself contains a full exceptional collection
of representations $(P_i)_{i \in Q_0}$, see e.g. \cite[Lemma 1.1]{MR1265279}.
The images under $\Phi_{\mathcal{U}}$ of these
$P_i$ are the summands $\mathcal{U}_{i}$ of the universal bundle of $\modulispace$, and form
a full exceptional collection for $\Phi_{\mathcal{U}}(\derived^{\bounded}(Q))$.

A positive answer to \cref{question:first-semiorthogonal-embedding}
yields thus an exceptional collection of $\#Q_0 + 1$ elements, and a positive answer to
\cref{question:r-semiorthogonal-embeddings} extends it to one with $r \cdot (\#Q_0 + 1)$
elements.
When this holds, it gives a finer decomposition of $\derived^{\bounded}(\modulispace)$
than the corresponding result in the case of curves.
Indeed, the derived category of a curve $\derived^{\bounded}(C)$ is indecomposable \cite[Theorem 1.1]{MR2838062},
and as a result
\eqref{equation:partial-SOD-VBAC} only gives a partial semiorthogonal decomposition of
$\derived^{\bounded}(\modulispace_{C}(r,\mathcal{L}))$.
In the case of quivers, not only do we give a partial semiorthogonal decomposition
of $\derived^{\bounded}(\modulispace)$ in \eqref{equation:partial-SOD},
but we further decompose each copy of $\derived^{\bounded}(Q)$ into $\#Q_0$ exceptional objects,
so that \cref{question:first-semiorthogonal-embedding,question:r-semiorthogonal-embeddings}
yield respectively the exceptional collections
\begin{equation}\label{equation:exceptional-collection-A}
    \mathcal{O}_{\modulispace}, \mathcal{U}_{1}, \dots, \mathcal{U}_{n}
\end{equation}
and
\begin{equation}\label{equation:exceptional-collection-B}
    \mathcal{O}_{\modulispace}, \mathcal{U}_{1}, \dots, \mathcal{U}_{n},
    \mathcal{O}(H), \mathcal{U}_{1}(H), \dots,\mathcal{U}_{n}(H), \ \ \dots, \ \
    \mathcal{O}((r-1)H), \mathcal{U}_{1}((r-1)H), \dots, \mathcal{U}_{n}((r-1)H).
\end{equation}
\end{remark}

\begin{remark}\label{remark:hochschild-homology}
In the case of curves, the semiorthogonal decomposition of \eqref{equation:partial-SOD-VBAC}
is known not to hold if $r = g = 2$. This can be proven using the additivity of
Hochschild homology: Kuznetsov \cite[Corollary 7.5]{kuznetsov2009hochschildhomologysemiorthogonaldecompositions}
has shown that, given a semiorthogonal decomposition
\begin{equation}
\derived^{\bounded}(X) = \langle \mathcal{C}_1,\dots, \mathcal{C}_n\rangle,
\end{equation}
the $i\mhyphen$th Hochschild homology of $X$ is given by
\begin{equation}\label{equation:Hochschild-homology-additive}
\Hochschild_{i}(X) = \bigoplus_{k = 1}^n \Hochschild_{i}(\mathcal{C}_k).
\end{equation}
We know that if $g = r = 2$, $\modulispace_{C}(2, \mathcal{L})$
is the intersection of two quadrics in $\PP^5$ (see \cite[Theorem~1]{MR237500}),
so $\dim \Hochschild_{0}(\modulispace_{C}(2, \mathcal{L})) = 4$;
we know that for a curve we always have
$\dim \Hochschild_{0}(C) = 2$,
and lastly we know that $\dim \Hochschild_{0}(\langle \Theta\rangle) = 1$, so
the decomposition~\eqref{equation:partial-SOD-VBAC} cannot be semiorthogonal if $r = g = 2$,
because the dimension of the right-hand side of~\eqref{equation:Hochschild-homology-additive}
in degree~$0$ would be~$6$.

In general, the dimension of $\Hochschild_{0}(X)$ gives the number of objects
that a full exceptional collection for~$\derived^{\bounded}(X)$, if it exists, must contain.
The Hochschild--Kostant--Rosenberg theorem \cite[Theorem 8.3]{kuznetsov2009hochschildhomologysemiorthogonaldecompositions}
allows to compute $\dim \Hochschild_{\bullet}(X)$ by knowing the Hodge diamond of $X$,
and for quiver moduli the latter is described in~\cite[Corollary~6.8]{MR1974891}
and implemented in \cite{hodge-diamond-cutter, quivertools}.
We can thus effectively compute these dimensions and detect at least some
cases where our questions must have a negative answer, similarly to the case $r = g = 2$
for curves.

In cases where \cref{question:first-semiorthogonal-embedding,question:r-semiorthogonal-embeddings}
have a positive answer, we also know how many objects our exceptional
sequence is missing to have a chance of being full.
\end{remark}

\begin{remark}\label{remark:frequency-pathological-behaviour}
While for vector bundles on curves the case $r = g = 2$ is the only known case where
\eqref{equation:partial-SOD-VBAC} does not hold (and it is in fact conjectured in
\cite{MR4651618} to be the \emph{only one}),
for quiver moduli several cases are known for which
\cref{question:first-semiorthogonal-embedding,question:r-semiorthogonal-embeddings}
have negative answers.
The evidence shown in \cref{section:good-cases} suggests that
examples with a negative answer stemming from \cref{remark:hochschild-homology} are
only to be expected in low dimension:
in such cases,
fewer objects than the amount predicted in
\cref{question:first-semiorthogonal-embedding,question:r-semiorthogonal-embeddings}
are required to form a full exceptional sequence.
Low-dimensional moduli of vector bundles on curves are scarce,
as $\dim(\modulispace_{C}(r,\mathcal{L})) = (r^2 - 1)(g - 1)$, whereas
quiver moduli in low dimensions are abundant:
under \cref{assumptions}, when fixing the canonical stability parameter~$\canonicalstability$,
$\modulispace^{\canonicalstability\stable}(Q,\tuple{d})$ is a smooth rigid Fano variety,
and already in dimensions $1$ and $2$ all rigid Fano varieties
(that is, $\PP^1$ and the six rigid del Pezzo surfaces) are realized.
The behaviour that is only observed once in
the case of curves is thus expected to appear more often for quiver moduli.
\end{remark}

In this paper, we will settle a weaker version of~\cref{question:r-semiorthogonal-embeddings}
for the canonical choice of stability parameter.
We will furthermore prove that
\cref{question:first-semiorthogonal-embedding,question:r-semiorthogonal-embeddings}
have a positive answer for a number of meaningful examples,
and illustrate known cases where they do not.
All the known cases where either question has a negative answer are predicted
by $\Hochschild_{0}(\modulispace)$ being too small, as explained
in \cref{remark:hochschild-homology}.


\paragraph{A conjecture of Schofield}

\begin{definition}
A vector bundle $\mathcal{E}$ in $\derived^{\bounded}(\modulispace)$ is called a \emph{tilting bundle}
if it classically generates $\derived^{\bounded}(\modulispace)$, and moreover
\begin{equation}\label{equation:partial-tilting-definition}
\Ext^{\geq 1}(\mathcal{E}, \mathcal{E}) = 0.
\end{equation}
When $\mathcal{E}$ satisfies \eqref{equation:partial-tilting-definition} but
does not necessarily generate the whole category, it is called a \emph{partial tilting bundle}.
\end{definition}

Given an exceptional collection of vector bundles $(A_1,\dots,A_n)$,
if no higher extensions exist among any two $A_i$ and $A_j$ the collection is said to
be \emph{strongly exceptional}; it is then immediate that $\bigoplus_{i=1}^{n} A_i$
is a partial tilting bundle.

We recall a result on tilting bundles on quiver moduli.
\begin{theorem}[Corollary C, \cite{rigidity-paper}]\label{theorem:partial-tilting}
Let $Q, \tuple{d}$ and $\theta$ satisfy \cref{assumptions}.
The universal representation $\mathcal{U}$
on $\modulispace^{\theta\stable}(Q,\tuple{d})$ is a partial tilting bundle, i.e.,
\begin{equation}
\Ext^{\geq 1}_{\modulispace^{\theta\stable}(Q,\tuple{d})}(\mathcal{U}, \mathcal{U}) = 0.
\end{equation}
\end{theorem}

This theorem partially settles a conjecture attributed to Schofield. A second
part of this conjecture~\cite[page~80]{MR1428456} states the following.

\begin{conjecture}\label{schofield-conjecture:complete-partial-tilting}
With notation as in \cref{theorem:partial-tilting}, the partial tilting bundle $\mathcal{U}$
can be completed to a tilting bundle.
\end{conjecture}

In \cref{section:good-cases}, we give a new proof of the well-known fact that,
for all the del Pezzo surfaces that can be realised as quiver moduli,
\cref{schofield-conjecture:complete-partial-tilting} holds.
In general, \cref{question:first-semiorthogonal-embedding,question:r-semiorthogonal-embeddings}
can be strengthened to attempt to give a partial answer to \cref{schofield-conjecture:complete-partial-tilting}:

\begin{alphaquestion}\label{question:strongly-exceptional}
Let $Q, \tuple{d}$ and $\theta$ be as in \cref{assumptions}.
Applying the decomposition of \eqref{equation:exceptional-collection-A}, respectively
\eqref{equation:exceptional-collection-B}, to
\eqref{equation:first-embedding}, respectively \eqref{equation:partial-SOD},
is the resulting exceptional collection strongly exceptional?
\end{alphaquestion}
A positive answer to \cref{question:strongly-exceptional} would give a partial tilting bundle for $\modulispace$
larger than the one given by \cref{theorem:partial-tilting}.
Indeed, the latter gives a bundle with $\#Q_0$ summands, while \cref{question:strongly-exceptional}
would yield larger ones, respectively with $\#Q_0 + 1$ or $r \cdot (\#Q_0 + 1)$ summands.

For weaker versions of \cref{question:r-semiorthogonal-embeddings,question:strongly-exceptional},
if one only considers the canonical stability parameter it is possible to give
a sufficient criterion to ensure an affirmative answer in general.
This technical criterion is related to the use of Teleman quantization,
and is stated in \cref{assumption:t-large}.

\begin{alphatheorem}\label{result:weaker-theorem}
Let $Q, \tuple{d}$ and $\canonicalstability$ satisfy \cref{assumptions}.
If \cref{assumption:t-large} is also satisfied,
there is a semiorthogonal decomposition of $\derived^{\bounded}(\modulispace)$ given by
\begin{equation}\label{equation:weaker-sod}
\left\langle
\Phi_{\mathcal{U}}(\derived^{\bounded}(Q)), \dots, \Phi_{\mathcal{U}((r-1)H)}(\derived^{\bounded}(Q)); \dots
\right\rangle.
\end{equation}
Furthermore, the exceptional collection obtained from it is strongly exceptional.
\end{alphatheorem}
\begin{remark}
\cref{result:weaker-theorem} is independent of the choice of linearisation defining $\mathcal{U}$.
\end{remark}
\begin{remark}
Considering this version of \cref{question:r-semiorthogonal-embeddings,question:strongly-exceptional}
has the advantage that the technical criterion in~\cref{assumption:t-large}
can be very efficiently verified using \cite{quivertools}.
By contrast, as discussed in \cref{section:good-cases}, verifying
\cref{question:first-semiorthogonal-embedding,question:r-semiorthogonal-embeddings,question:strongly-exceptional}
for a given example involves computations in the Chow ring,
for which the complexity can become prohibitively high.
\end{remark}

In relevant cases, such as for $m$-Kronecker quivers with dimension vector $(2, 3)$,
the technical criterion of \cref{assumption:t-large} can be verified directly for any
$m \in \ZZ_{\geq 3}$. This is in fact the proof of the first part of \cref{proposition:m-kronecker-quivers}.

The structure of this paper is the following:
in \cref{section:quiver-moduli} we recall the standard construction of quiver moduli.
\cref{section:embedding-multiple-copies} describes the functor $\Phi_{\mathcal{U}}$
and formulates an equivalent criterion for \eqref{equation:first-embedding} and
\eqref{equation:partial-SOD} to be exceptional collections. \cref{section:HN-Teleman}
explains Teleman quantization, one of the two main ingredients in our proofs.
Lastly, the necessary intersection theory and our results
are presented in \cref{section:good-cases}.

The computations that we perform in the Chow ring are made possible by the
effective description of Chow rings of quiver moduli in \cite{MR3318266, MR4770368}.
Teleman quantization is a theorem about equivariant cohomology
of $G$-linearised quasicoherent sheaves on $G$-varieties.
We state a version of it adapted to our setting in \cref{theorem:teleman-quantization},
which was first formulated for the quiver moduli setting in \cite{rigidity-paper}.

\subsection*{Acknowledgements}
We wish to thank Pieter Belmans for suggesting
\cref{question:first-semiorthogonal-embedding,question:r-semiorthogonal-embeddings,question:strongly-exceptional}
and for his many ideas and discussions about this work.
We wish to thank the referees for their valuable suggestions and comments.
Preliminary results on these topics were first presented at the conference
\href{https://www.mis.mpg.de/de/events/event/decomposing-quiver-moduli-a-quivertools-showcase}{MEGA 2024},
at the Max Plank Institute for Mathematics in the Sciences,
in Leipzig, Germany, on July 29th, 2024.

\paragraph{Funding} This work was supported by the Luxembourg National Research Fund (FNR-17953441).

\section{Moduli spaces of quiver representations}\label{section:quiver-moduli}
A \emph{quiver} $Q$ is a finite directed graph composed of vertices $Q_0$ and arrows $Q_1$.
Each arrow $a \in Q_1$ has a \emph{source} $s(a) \in Q_0$ and
a \emph{target} $t(a) \in Q_0$.
Vertices in $Q_0$ are numbered from $1$ to $n$.
We work over an algebraically closed field $\field$ of characteristic $0$.

A \emph{representation} of the quiver $Q$ is the data of a finite-dimensional vector space
$V_i$ for each vertex $i \in Q_0$ and a linear morphism
$V_a \colon V_{s(a)} \to V_{t(a)}$ for each arrow $a \in Q_1$.
The dimensions of the $V_i$ are encoded in a \emph{dimension vector}
$\tuple{d} = (d_1,\dots,d_n) \in \NN^{Q_0}$, and
each representation is identified with a point in the \emph{parameter space}
\begin{equation}
\parameterspace\colonequals
\mathrm{R}(Q,\tuple{d})\colonequals
\bigoplus_{a \in Q_1} \mathrm{Mat}_{d_{t(a)}\times d_{s(a)}}(\field).
\end{equation}
The group $\GL_{\tuple{d}}\colonequals \prod_{i \in Q_0} \GL_{d_i}(\field)$
acts on $\parameterspace$ by conjugation, and its orbits correspond to
isomorphism classes of representations.
Once we fix a \emph{stability parameter} $\theta \in \ZZ^{Q_0} \setminus{\{\tuple{0}\}}$
such that $\theta \cdot \tuple{d} = 0$, we define a representation to be \emph{stable},
respectively \emph{semistable}, if all its nontrivial subrepresentations $W$ satisfy
the inequality $\theta \cdot \dim(W) < 0$, respectively $\theta \cdot \dim(W) \leq 0$.

We can take the GIT quotient of the $\theta$-stable locus
$\parameterspace^{\theta\stable}(Q,\tuple{d}) \geomquotient_{\theta}\, \GL_{\tuple{d}}$;
denoted by $ \modulispace\colonequals\modulispace^{\theta\stable}(Q,\tuple{d})$,
this quotient is the \emph{moduli space of $\theta$-stable representations}
of $Q$ of dimension vector $\tuple{d}$.

Given $Q$ and a dimension vector $\tuple{d}$, one defines the
\emph{canonical stability parameter} $\canonicalstability$
as the linear function sending $\tuple{e}$ to
$\langle\tuple{d}, \tuple{e}\rangle - \langle \tuple{e}, \tuple{d} \rangle$.
Here the brackets denote the bilinear Euler form of the quiver.

Given a tuple of integers $\tuple{e} \in \ZZ^{Q_0}$
such that $\tuple{e} \cdot \tuple{d} = 0$,
we can define a linearisation of the trivial
line bundle $\mathcal{O}_{\parameterspace^{\theta\semistable}}(Q, \tuple{d})$ that
descends to the GIT quotient $\modulispace^{\theta\semistable}(Q, \tuple{d})$.
We denote the linearised line bundle and its descent by $L(\tuple{e})$
and $\mathcal{L}(\tuple{e})$, respectively.

We introduce several assumptions which will be used throughout this work.
\begin{assumption}\label{assumptions} \phantom{aa}
    \begin{enumerate}
        \item
        The quiver $Q$ is acyclic,
        \item
        For all subdimension vectors $\tuple{d'} \leq \tuple{d}$, we have $\tuple{d'}\cdot\theta \neq 0$, and
        \item For all subdimension vectors $\tuple{d'} \leq \tuple{d}$ satisfying $\theta \cdot \tuple{d'} > 0$,
        the bilinear Euler form $\langle -,-\rangle$ yields $\langle \tuple{d'}, \tuple{d - d'}\rangle \leq -2$.
    \end{enumerate}
The last technical condition is referred to as \emph{$\theta$-strong ample stability},
see \cite[Definition 4.1]{rigidity-paper}.
It was first introduced in \cite[Proposition 5.1]{MR3683503},
where it is used as a sufficient condition
to establish that the codimension of the unstable locus
$\parameterspace(Q,\tuple{d})\setminus\parameterspace^{\theta\semistable}(Q,\tuple{d})$
is at least $2$.
\end{assumption}

\begin{remark}
The second condition is known as $\theta$-coprimality.
Assuming that a dimension vector $\tuple{d}$ is $\theta$-coprime
implies that $\theta$-semistable representations are actually $\theta$-stable,
so that stability and semistability are equivalent.
\end{remark}

Under \cref{assumptions}, the moduli space is a smooth projective variety which comes
with a universal bundle
\begin{equation}\label{equation:universal-bundles}
\mathcal{U} = \bigoplus_{i \in Q_0} \mathcal{U}_i,
\end{equation}
such that for each point $[V] \in \modulispace$ corresponding to
the isomorphism class of the representation $V$, the fiber $\mathcal{U}_{[V]}$
is a representation isomorphic to $V$ itself.

\begin{remark}
When a quiver is acyclic we fix a numbering of the vertices to give a total order on them.
\end{remark}

\begin{remark}\label{remark:choice-linearisation}
The summands of $\mathcal{U}$ are not unique: they depend on a choice of a linearisation,
that is, a tuple $\tuple{a} \in \ZZ^{Q_0}$ for which $ \tuple{a} \cdot \tuple{d} = 1$.
When $\tuple{a}$ is not omitted, the notation $\mathcal{U}_i(\tuple{a})$ is used.
Given two such choices $\tuple{a}$ and $\tuple{a'}$,
the resulting bundles $\mathcal{U}_{i}(\tuple{a})$ and $\mathcal{U}_{i}(\tuple{a'})$
are isomorphic up to tensoring by an equivariant line bundle:
\begin{equation}
\mathcal{U}_{i}(\tuple{a}) \simeq \mathcal{U}_{i}(\tuple{a'})
\otimes \mathcal{L}(\tuple{a - a'}).
\end{equation}
\end{remark}

We set up an example which will be used through the paper
to illustrate various notions.

\begin{example}\label{example:recurring-example}
Let $Q$ be the 3-Kronecker quiver, $\tuple{d} = (3, 4)$ and $\theta = (12,-9)$.
\begin{equation}
    \begin{tikzpicture}[baseline = -3pt, node distance=1.5cm]
        \node (1)                   {$\bullet$};
        \node (2) [right of = 1]    {$\bullet$};

        \draw (1) node[above] {1};
        \draw (2) node[above] {2};

        \draw[->] (1) edge (2);
        \draw[->, bend  left = 30] (1) edge (2);
        \draw[->, bend right = 30] (1) edge (2);
    \end{tikzpicture}
\end{equation}
It is straightforward to verify that \cref{assumptions} holds.
This can also be done using \textsc{QuiverTools}.

The resulting moduli space $\modulispace^{\theta\stable}(Q, \tuple{d})$ is a smooth
projective 12-dimensional rational Fano variety~\cite[Theorem 4.2]{MR4352662}.
Its Hodge diamond only has nonzero entries on the central column:
\begin{equation}\label{equation:hodge-column-recurring-example}
(\hh^{0,0},\dots,\hh^{12, 12}) = (1, 1, 3, 5, 8, 10, 12, 10, 8, 5, 3, 1, 1).
\end{equation}
The expected length of a full exceptional collection for
$\modulispace^{\theta\stable}(Q, \tuple{d})$ is thus $\sum_{i = 0}^{12} \hh^{i,i} = 68$.
\end{example}

\section{Semiorthogonal embeddings}
\label{section:embedding-multiple-copies}
In this section we describe the fully faithful functors that embed $\derived^{\bounded}(Q)$
into $\derived^{\bounded}(\modulispace)$, and give equivalent criteria for
\cref{question:first-semiorthogonal-embedding,question:r-semiorthogonal-embeddings,question:strongly-exceptional}
to have positive answers.

The object $\mathcal{U}$ is a representation of $Q$ in the
category of quasicoherent sheaves on $\modulispace$, that is,
it admits a structure of right $kQ$-module and
a structure of left $\mathcal{O}_{\modulispace}$-module.
We define the Fourier--Mukai functor
\begin{equation}\label{equation:definition-FMT}
    \Phi_{\mathcal{U}} : \derived^\bounded(Q) \to
    \derived^\bounded(\modulispace) :
    V \mapsto \mathcal{U} \otimes^\mathrm{L}_{kQ} V.
\end{equation}
For any line bundle $\mathcal{O}(D)$ on $\modulispace$,
we also define the twisted Fourier--Mukai functor
\begin{equation}\label{equation:twisted-FMT}
\Phi_{\mathcal{U}(D)} \colonequals \left(-\otimes \,\mathcal{O}(D)\right)\,\circ\, \Phi_{\mathcal{U}}.
\end{equation}

As twisting by a line bundle is an autoequivalence of categories,
we state the following slight generalization of a result of \cite{vector-fields-paper}.
\begin{lemma}[Theorem D, \cite{vector-fields-paper}]\label{lemma:multiple-embeddings}
Let $Q$, $\tuple{d}$ and $\theta$ as in \cref{assumptions},
and let $D$ be a divisor on $\modulispace$.
The functor $\Phi_{\mathcal{U}(D)}$ is fully faithful,
and thus embeds a copy of $\derived^\bounded(Q)$ in $\derived^\bounded(\modulispace)$.
\end{lemma}

Denote by $r$ the Mukai index of $\modulispace$,
i.e., the largest integer $s$ such that $K_{\modulispace} = s\cdot L$ for some divisor $L$,
and let $H \colonequals - \frac{1}{r}K_{\modulispace}$ be the irreducible part of the
anticanonical divisor.
In this work we will be concerned with twisting the Fourier--Mukai functor
above by multiples of $\mathcal{O}(H)$.

In the rest of this section we establish equivalent criteria for
\cref{question:first-semiorthogonal-embedding,question:r-semiorthogonal-embeddings,question:strongly-exceptional}
to have positive answers.

A semiorthogonal decomposition of $\derived^{\bounded}(\modulispace)$
is an ordered sequence of full triangulated admissible subcategories that are,
in the given order, left orthogonal.
More precisely, if there are $n$ full triangulated admissible subcategories
$\mathcal{A}_i \subset \derived^{\bounded}(\modulispace)$ such that
\begin{equation}\label{def:SOD-definition}
\derived^{\bounded}(\modulispace) = \left\langle \mathcal{A}_1,\dots,\mathcal{A}_n \right\rangle,
\end{equation}
then \eqref{def:SOD-definition} is said to be a \emph{semiorthogonal decomposition} if
for all $E_i \in \mathcal{A}_i$ and all $E_j \in \mathcal{A}_j$
with $i < j$ and for all $k \in \ZZ$ one has
\begin{equation}\label{equation:define-sod}
\Hom(E_j, E_i[k]) = 0.
\end{equation}

For \cref{question:r-semiorthogonal-embeddings} to hold it
is then necessary and sufficient that for any $0 \leq n_1 < n_2 \leq r-1$,
for any $V, W \in \derived^\bounded(Q)$ and for any $k \in \ZZ$,
we have
\begin{align}
    \Hom(\Phi_{\mathcal{U}(n_2H)}(V),\Phi_{\mathcal{U}(n_1H)}(W)[k]) &= 0, \label{SO-condition:Db-Db}\\
    \Hom(\mathcal{O}(n_2H), \mathcal{O}(n_1H)[k]) &= 0, \label{SO-conditions:O-O} \\
    \Hom(\mathcal{O}(n_2H), \Phi_{\mathcal{U}(n_1H)}(V)[k]) &= 0, \label{SO-conditions:O-Db}
\end{align}
and that for any $0 \leq n_1 \leq n_2 \leq r-1$ we have
\begin{equation}
    \Hom(\Phi_{\mathcal{U}(n_2H)}(V), \mathcal{O}(n_1H)[k]) = 0. \label{SO-conditions:Db-O}
\end{equation}
We give a sufficient condition to ensure that these vanishings are satisfied.

\begin{lemma}\label{lemma:criterion-QB}
Let $Q, \tuple{d}$ and $\theta$ be as in \cref{assumptions}.
If we can find a linearisation $\tuple{a}$ for $\mathcal{U}$ such that, for all
$k \geq 0$, for all $1 \leq s \leq r-1$ and for all $i,j \in Q_0$,
\begin{align}
    \mathrm{H}^k(\modulispace, \mathcal{U}^\vee_i \otimes \mathcal{U}_j \otimes \mathcal{O}(-sH)) &= 0, \label{equation:needed-UidualUjH}\\
    \mathrm{H}^k(\modulispace, \mathcal{O}(-sH)) &= 0, \label{equation:needed-H}\\
    \mathrm{H}^k(\modulispace, \mathcal{U}_i \otimes \mathcal{O}(-sH)) &= 0, \label{equation:needed-UiH}
\end{align}
and that for all $0 \leq t \leq r-1$,
\begin{equation}
    \mathrm{H}^k(\modulispace, \mathcal{U}^\vee_i \otimes \mathcal{O}(-tH)) = 0, \label{equation:needed-UidualH}
\end{equation}
then \cref{question:r-semiorthogonal-embeddings} has a positive answer.
\end{lemma}
\begin{proof}
By \cite[Lemma 25]{MR3954042}, to check conditions
\eqref{SO-condition:Db-Db}, \eqref{SO-conditions:Db-O}
and \eqref{SO-conditions:O-Db} it is enough to
verify the vanishings for the elements of a
spanning class of $\derived^\bounded(Q)$.
Using the fact that the indecomposable projective representations
$\{P_i\}_{i \in Q_0}$ of $Q$ form such a family,
and that the image of $P_i$ via $\Phi_{\mathcal{U}}$ is $\mathcal{U}_{i}$,
we conclude.
\end{proof}

\begin{corollary}\label{corollary:criterion:QA}
If \eqref{equation:needed-UidualH} holds for all $i$, for all $k$ and for $t = 0$,
then \cref{question:first-semiorthogonal-embedding} has a positive answer.
\end{corollary}

The same idea can be used to check whether an exceptional collection
is strongly exceptional, as we state below.
\begin{lemma}\label{lemma:criterion:QC}
Under the assumptions of \cref{lemma:criterion-QB}, if for all
$1 \leq s \leq r-1$ and for all $i,j \in Q_0$,
\begin{align}
    \mathrm{H}^{\geq 1}(\modulispace, \mathcal{U}^\vee_i \otimes \mathcal{U}_j \otimes \mathcal{O}(sH)) &= 0,\\
    \mathrm{H}^{\geq 1}(\modulispace, \mathcal{O}(sH)) &= 0, \\
    \mathrm{H}^{\geq 1}(\modulispace, \mathcal{U}^{\vee}_i \otimes \mathcal{O}(sH)) &= 0,
\end{align}
and that for all $0 \leq t \leq r-1$,
\begin{equation}
    \mathrm{H}^{\geq 1}(\modulispace, \mathcal{U}_i \otimes \mathcal{O}(tH)) = 0,
\end{equation}
then the collection of \cref{question:strongly-exceptional} is strongly exceptional.
\end{lemma}

\begin{remark}
The terms in \eqref{equation:needed-UidualUjH} and
\eqref{equation:needed-H} do not depend on the choice of linearisation explained in
\cref{remark:choice-linearisation}, while \eqref{equation:needed-UiH} and
\eqref{equation:needed-UidualH} do.
\end{remark}

To show these vanishings in many relevant examples, we will use a combination of
three ingredients: Teleman quantization, Chow ring computations and Serre duality.
We establish the following lemma for future reference.

\begin{lemma}\label{lemma:serre-duality}
Let $Q, \tuple{d}$ and $\theta$ be as in \cref{assumptions}. Denote as before by $r$ the index
of $\modulispace$. For all $i, j \in Q_0$, for all $k \in \ZZ$ and for all $s = 0,\dots, r$
we have the isomorphisms
\begin{align}
\HH^{k}(\modulispace, \mathcal{U}^{\vee}_i \otimes \mathcal{U}_j \otimes \mathcal{O}(-sH)) &\simeq
\HH^{\dim(M) - k}(\modulispace, \mathcal{U}^{\vee}_j \otimes \mathcal{U}_i \otimes \mathcal{O}(-(r - s)H))^{\vee}, \\
\HH^{k}(\modulispace, \mathcal{U}^{\vee}_i \otimes \mathcal{O}(-sH)) &\simeq
\HH^{\dim(M) - k}(\modulispace, \mathcal{U}_i \otimes \mathcal{O}(-(r-s)H))^{\vee}.
\end{align}
\end{lemma}
\begin{proof}
This is an immediate application of Serre duality, considering that
$K_{\modulispace} = -rH$ by definition.
\end{proof}

\section{Stratification of the unstable locus and Teleman quantization}\label{section:HN-Teleman}

In this section we describe a stratification of the $\theta$-unstable locus of
$\parameterspace(Q, \tuple{d})$ and state a version of Teleman quantization
adapted to the scope of this work.
For our purposes, Teleman quantization will be used as a sufficient criterion
to establish vanishing of higher cohomology for vector bundles on quiver moduli.
It is in fact a deep result of Teleman \cite{MR1792291}, later extended to the setting
of derived categories by Halpern-Leistner~\cite{MR3327537}.

Fix a quiver $Q$ and a stability parameter $\theta$.
Let \emph{the slope function} $\mu_{\theta}$ be defined by
\begin{equation}
    \mu_{\theta} \colon \ZZ^{Q_0} \to \QQ: \tuple{e} \mapsto \frac{\theta \cdot \tuple{e}}{\Sigma_{i \in Q_0} e_i}.
\end{equation}

We say that a representation $V$ of $Q$ is $\mu_{\theta}$-stable, respectively semistable,
if all its nonzero proper subrepresentations $W$ satisfy the inequality $\mu_{\theta}(\dim(W)) < 0$,
respectively $\mu_{\theta}(\dim(W)) \leq 0$.
Note that we do not require that $\mu_{\theta}(\dim(V)) = 0$.

\begin{definition}\label{definition:HN-type}
Let $Q, \tuple{d}$ and $\theta$ be a quiver, a dimension vector and a stability parameter.
A \emph{Harder--Narasimhan type} $\tuple{d}^*$ is an ordered sequence of dimension
vectors $(\tuple{d}^1,\dots,\tuple{d}^{\ell})$ that sum to $\tuple{d}$, such that
$\mu(\tuple{d}^m) > \mu(\tuple{d}^{m+1})$ for all $m = 1,\dots,\ell - 1$, and such that for each $m = 1, \dots, \ell$
there exists a $\mu_{\theta}$-semistable representation of $Q$ of dimension vector $\tuple{d}^{m}$.
\end{definition}

It is a theorem of Reineke \cite[Proposition 2.5]{MR1974891} that any
representation $V \in \parameterspace(Q,\tuple{d})$ admits a unique filtration
$0 = V_0 \subset V_1 \dots \subset V_{\ell} = V$ such that $(\dim(V_{m + 1}/V_{m}))_{m = 0}^{\ell-1}$
is a Harder--Narasimhan type.
This allows to say that each representation in $\parameterspace(Q, \tuple{d})$
\emph{has a unique Harder--Narasimhan type}.
From the definition it is clear that, for a stability parameter $\theta$,
a representation is $\theta$-semistable if and only if its Harder--Narasimhan type is
$(\tuple{d})$, i.e., it is trivial.

It is another theorem of Reineke \cite{MR1974891} that
the $\theta$-\emph{unstable locus} - that is, the complement
of~$\parameterspace^{\theta\semistable}(Q, \tuple{d})$,
admits a stratification into locally closed, smooth, disjoint subsets $S_{\tuple{d}^*}$,
indexed by the Harder--Narasimhan types $\tuple{d}^*$ and characterized by the fact
that all the representations in $S_{\tuple{d}^*}$  have Harder--Narasimhan type
$\tuple{d}^*$.

\begin{example}
In the case of \cref{example:recurring-example}, there are 19 Harder--Narasimhan types; one is the
trivial one $((3, 4))$, corresponding to the semistable locus, and the 18 others correspond
to the strata of the unstable locus. The nontrivial Harder--Narasimhan types are
listed in the first column of \cref{tab:HN-recurring-example}.
\end{example}

To each nontrivial Harder--Narasimhan type $\tuple{d}^* = (\tuple{d}^1,\dots,\tuple{d}^{\ell})$
we associate a one-parameter subgroup of $\GL_{\tuple{d}}$ as follows:
let $k_m \colonequals c\cdot\mu(\tuple{d}^m)$,
where $c$ is the smallest natural number for which, for all $m$, $k_m \in \ZZ$.
We define the one-parameter subgroup $\lambda_{\tuple{d}^*}$ to be,
on each component $i \in Q_0$, the diagonal matrix
\begin{equation}\label{equation:1-PS}
\lambda_{\tuple{d^*}, i}(z) \colonequals \operatorname{diag}\big(
        \underbrace{z^{k_1},\ldots,z^{k_1}}_{d_i^1 \text{ times}},
        \underbrace{z^{k_2},\ldots,z^{k_2}}_{d_i^2 \text{ times}},
        \ldots,
        \underbrace{z^{k_\ell},\ldots,z^{k_\ell}}_{d_i^\ell \text{ times}}
    \big).
\end{equation}

Let now $\{S_{\tuple{d}^*}\}$ be the Harder--Narasimhan
stratification of the $\theta$-unstable locus, and for all types $\tuple{d}^*$
let $\lambda_{\tuple{d}^*}$ be defined as above.
Let $Z_{\tuple{d}^*}$ be the set of points in $S_{\tuple{d}^*}$
fixed by $\lambda_{\tuple{d}^*}$, and
denote the weight of the natural action of $\lambda_{\tuple{d}^*}$ on
$\det(\mathrm{N}^\vee_{S_{\tuple{d}^*}/R})|_{Z_{\tuple{d}^*}}$
by $\eta_{\tuple{d}^*}$.
Consider a $\GL_{\tuple{d}}$-linearised quasicoherent sheaf on $\parameterspace$, say $F$,
and assume that it descends to a quasicoherent sheaf $\mathcal{F}$ on the quotient $\modulispace$.
On each $Z_{\tuple{d}^*}$, denote the set of weights of $\lambda_{\tuple{d}^*}$ on~$F|_{Z_{\tuple{d}^*}}$ by $W(F, \tuple{d}^*)$.

We state the Teleman quantization theorem in a way adapted to the scope of this work,
rather than in its most general form.
\begin{theorem}[Theorem 3.6 \cite{rigidity-paper}]\label{theorem:teleman-quantization}
    If the inequality
    \begin{equation}\label{equation:teleman-inequality}
        \max(W(F,\tuple{d}^*)) < \eta_{\tuple{d}^*}
    \end{equation}
    holds for every Harder--Narasimhan type $\tuple{d}^*$, then for all $k \geq 0$ there is
    an isomorphism
\begin{equation}\label{equation:teleman-theorem}
    \HH^k(\modulispace, \mathcal{F}) \cong \HH^k(\parameterspace, F)^{\GL_{\tuple{d}}}.
\end{equation}
\end{theorem}
In \eqref{equation:teleman-theorem}, the right-hand side is the $\GL_{\tuple{d}}$-invariant
subgroup of $\HH^k(\parameterspace, F)$. In particular, since $\parameterspace$ is affine,
we obtain the vanishing result below.
\begin{corollary}\label{corollary:teleman-vanishing}
If \eqref{equation:teleman-inequality} holds on all Harder--Narasimhan strata,
for all $k \geq 1$ we have $\HH^k(\modulispace, \mathcal{F}) = 0$.
\end{corollary}

Below we collect the weights for the vector bundles needed for this work.

\begin{proposition}[Lemma 3.19 \cite{rigidity-paper}]
    \label{proposition:weights-universal-bundles}
Let $Q, \tuple{d}$ and $\theta$ be as in \cref{assumptions}.
Let $\tuple{a}$ be the linearisation defining the bundles $U_i(\tuple{a})$ on $\parameterspace$.
Then, on the stratum $S_{\tuple{d}^*}$, the $\lambda_{\tuple{d}^*}$-weights of $U_i(a)$ are given by
\begin{equation}
    W(U_i(\tuple{a}),\tuple{d}^*) = \left\{k_m - \sum_{j \in Q_0}\sum_{n = 1}^{\ell}a_j d^n_j k_n ~|~ 1 \leq m \leq \ell\right\},
\end{equation}
where the weight indexed by $m$ appears with multiplicity $d_i^m$.
\end{proposition}

An immediate corollary is that the weights of $U^{\vee}_i \otimes U_j$ do not depend
on a choice of linearisation $\tuple{a}$.
\begin{corollary}\label{corollary:weight-uidualuj}
The $\lambda_{\tuple{d}^*}$-weights of $U^{\vee}_i \otimes U_j$
on the stratum $S_{\tuple{d}^*}$ are given by
\begin{equation}\label{corollary:weights-endomorphism-bundles}
    W(U^{\vee}_i \otimes U_j,\tuple{d}^*) = \left\{ k_m - k_n ~|~ 0 \leq m,n \leq \ell \right\},
\end{equation}
where the weight indexed by $m,n$ appears with multiplicity $d_{j}^{m} \cdot d_{i}^{n}$.
\end{corollary}

To describe $H$ and compute its weights, we use \cite[Proposition 4.1]{MR4352662},
which under \cref{assumptions} describes explicitly the linearised line bundle $L(\canonicalstability)$ on
$\parameterspace$ from which the canonical divisor descends.
From this, we know that~$H$ descends from
$L(\frac{1}{r}\canonicalstability)$.
In particular, $r = \gcd(\canonicalstability)$.

\begin{proposition}\label{proposition:weight-H}
In the same setting as \cref{proposition:weights-universal-bundles},
on the stratum $S_{\tuple{d}^*}$, the $\lambda_{\tuple{d}^*}$-weight of
the line bundle $L(\frac{1}{r}\canonicalstability)$ is
\begin{align}
   W(L(\tfrac{1}{r}\canonicalstability), \tuple{d}^*)
   &= \left\{-\frac{1}{r} \sum_{i \in Q_0} \canonicalstability_{i}~\sum_{n = 1}^\ell d^n_i k_n\right\} \\
   &= \left\{-\frac{1}{r} \canonicalstability \cdot \sum_{n = 1}^\ell \tuple{d}^n k_n\right\}.
\end{align}
\end{proposition}

\begin{remark}\label{remark:H-negative-weights}
As $k_s = \frac{\canonicalstability \cdot \tuple{d}^s}{\sum_i d^s_i}$,
the Teleman weights of $H$ are all strictly negative.
\end{remark}

Lastly, the weights $\eta_{\tuple{d}^*}$ are given below.

\begin{proposition}[Corollary 3.18 \cite{rigidity-paper}]\label{proposition:Teleman-bounds}
In the same setting as \cref{proposition:weights-universal-bundles},
on the stratum $S_{\tuple{d}^*}$, the \emph{Teleman bound} $\eta_{\tuple{d}^*}$ is
\begin{equation}
\eta_{\tuple{d}^*} = \sum_{1 \leq s<t\leq \ell}(k_{t} - k_{s})\langle \tuple{d}^s, \tuple{d}^t \rangle.
\end{equation}
\end{proposition}

\begin{remark}
As the slopes of the terms in $\tuple{d}^*$ are strictly decreasing and
each term admits $\mu$-semistable representations,
we have~$\hom(\tuple{d}^s,\tuple{d}^t) = 0$ for all $s < t$.
Therefore, $\langle \tuple{d}^s, \tuple{d}^t \rangle = -\ext(\tuple{d}^s, \tuple{d}^t)$.
\end{remark}

\begin{example}
For \cref{example:recurring-example}, we fix the linearisation $\tuple{a} = (3, -2)$
and use methods from \textsc{QuiverTools} \cite{quivertools} to compute relevant weights
and collect them in \cref{tab:HN-recurring-example}.

\begin{table}[H]
\centering
\begin{tabular}{ccccccc}
\toprule
$\tuple{d}^*$ &
$W(U^{\vee}_i \otimes U_j, \tuple{d}^*)$ &
$W(U_i(\tuple{a}), \tuple{d}^*)$ &
$W(L(\frac{\canonicalstability}{r}), \tuple{d}^*)$ &
$\eta_{\tuple{d}^*}$ \\
\midrule
((1, 1), (2, 3))                 & -21 0 21                           & -21 0                 & -21    & 84    \\
((2, 2), (1, 2))                 & -7 0 7                             & -14 -7                & -14    & 42    \\
((3, 3), (0, 1))                 & -21 0 21                           & -63 -42               & -63    & 126   \\
((3, 2), (0, 2))                 & -63 0 63                           & -315 -252             & -378   & 882   \\
((2, 1), (1, 3))                 & -35 0 35                           & -140 -105             & -175   & 455   \\
((2, 1), (1, 2), (0, 1))         & -14 -7 0 7 14                      & -49 -42 -35           & -56    & 133   \\
((2, 1), (1, 1), (0, 2))         & -28 -21 -7 0 7 28                  & -133 -112 -105        & -161   & 385   \\
((3, 1), (0, 3))                 & -63 0 63                           & -441 -378             & -567   & 1512  \\
((1, 0), (2, 4))                 & -14 0 14                           & -42 -28               & -56    & 140   \\
((1, 0), (2, 3), (0, 1))         & -105 -63 -42 0 42 63 105           & -315 -273 -210        & -378   & 882   \\
((1, 0), (1, 1), (1, 3))         & -63 -42 -21 0 21 42 63             & -210 -189 -147        & -273   & 693   \\
((1,\,0),\,(1,\,1),\,(1,\,2),\,(0,\,1)) & -42 -28 -21 -14 -7 0 7 14 21 28 42 & -133 -119 -112 -91    & -161   & 385   \\
((1, 0), (2, 2), (0, 2))         & -42 -21 0 21 42                    & -168 -147 -126        & -210   & 504   \\
((1, 0), (2, 1), (0, 3))         & -21 -14 -7 0 7 14 21               & -119 -105 -98         & -154   & 406   \\
((2, 0), (1, 3), (0, 1))         & -84 -63 -21 0 21 63 84             & -441 -420 -357        & -567   & 1512  \\
((2, 0), (1, 2), (0, 2))         & -21 -14 -7 0 14 21                 & -119 -112 -98         & -154   & 406   \\
((2, 0), (1, 1), (0, 3))         & -42 -21 0 21 42                    & -273 -252 -231        & -357   & 966   \\
((3, 0), (0, 4))                 & -21 0 21                           & -189 -168             & -252   & 756   \\
\bottomrule
\end{tabular}
\caption{Harder--Narasimhan types and Teleman weights for \cref{example:recurring-example}.}
\label{tab:HN-recurring-example}
\end{table}
\end{example}

\section{Computations and proofs of our claims}\label{section:good-cases}

In the following section we will give answers to
\cref{question:first-semiorthogonal-embedding,question:r-semiorthogonal-embeddings,question:strongly-exceptional}
for various relevant examples.

The intersection theory on quiver moduli is described in \cite{MR3318266,MR4770368},
which give a presentation of the Chow ring of $\modulispace^{\theta\stable}(Q, \tuple{d})$,
compute the point class and the Todd class. This gives an effective way to apply
the Hirzebruch--Riemann--Roch theorem and compute Euler characteristics.
Below we briefly collect the notions and results necessary
to compute Euler characteristics for the terms appearing in \cref{lemma:criterion-QB}.

\begin{theorem}[Franzen, \cite{MR3318266}]
Let $Q$ be an acyclic quiver, $\tuple{d}$ a dimension vector and $\theta$ a stability
parameter for which $\tuple{d}$ is $\theta$-coprime.
The Chow ring of the moduli space $\modulispace \colonequals \modulispace^{\theta\stable}(Q, \tuple{d})$ is given by
\begin{equation}
\mathrm{CH}^{\bullet}(\modulispace) = \QQ[\xi_{i, k_i}]_{i \in Q_0, k_i = 1,\dots,d_i} / I,
\end{equation}
where $\xi_{i, j} = \textrm{c}_{j}(\mathcal{U}_i)$, and $I$ is defined by several tautological
relations and by the linear relation $ \sum_{i}\xi_{i, 1} a_i = 0$.
\end{theorem}

\paragraph{Note}
The linear relation in $I$ is not tautological, as it depends on the choice
of linearisation $\tuple{a}$ described in \cref{remark:choice-linearisation}.
For more details the reader is directed to \cite[Section 2.2]{MR4770368}.

Under \cref{assumptions}, the first Chern class of the canonical bundle
is described in \cite[Remark 4.1]{MR4352662}
in terms of the Euler form of the quiver:
\begin{equation}
c_1(\omega_{\modulispace}) =
\sum_{i \in Q_0}
\left(\langle \tuple{d}, \tuple{i}\rangle - \langle \tuple{i},\tuple{d} \rangle\right)
\xi_{i, 1}.
\end{equation}
Using this, we can compute the Chern character of $\mathcal{O}(H)$.

\subsection{\texorpdfstring{$m$}{m}-Kronecker quivers}
Let $Q$ be the $m$-Kronecker quiver. For any $\tuple{d} = (d_1, d_2)$,
we have $\canonicalstability = m \cdot (d_2, -d_1)$, so if $\tuple{d}$ is $\canonicalstability$-coprime,
the index of~$\modulispace^{\theta\stable}(Q, \tuple{d})$ equals~$m$.
\begin{equation}
\label{equation:m-Kronecker-quiver}
    \begin{tikzpicture}[node distance=0.8cm]
        \node (1)                   {$\bullet$};
        \node (2) [right of = 1]    {$\vdots$};
        \node (3) [right of = 2]    {$\bullet$};

        \draw[->, bend  left = 20] (1) edge (3);
        \draw[->, bend  left = 30] (1) edge (3);
        \draw[->, bend right = 33] (1) edge (3);
    \end{tikzpicture}
    \end{equation}

We show that if we fix the dimension vector to be $\tuple{d} = (2, 3)$,
the decomposition in \eqref{equation:partial-SOD} always gives a partial tilting bundle
on $\modulispace^{\canonicalstability\stable}(Q,\tuple{d})$.

\begin{proposition}\label{proposition:m-kronecker-quivers}
Let $Q$ be a $m$-Kronecker quiver. Consider the dimension
vector~$\tuple{d} = (2, 3)$, the stability parameter $\canonicalstability$ and
the choice of linearisation $\tuple{a} = (2, -1)$.
The moduli space $\modulispace \colonequals \modulispace^{\theta\stable}(Q, \tuple{d})$ admits a
strongly exceptional collection given by
\begin{equation}\label{equation:exceptional-collection-2-3-m-kronecker}
\mathcal{U}_1, \mathcal{U}_2,
\mathcal{U}_1(H), \mathcal{U}_2(H), \dots,
\mathcal{U}_1((m-1)H), \mathcal{U}_2((m-1)H).
\end{equation}
If, moreover, we assume that for all $i \in Q_0$, $\HH^{0}(\modulispace, \mathcal{U}^{\vee}_{i}) = 0$,
a larger strongly exceptional collection is given by
\begin{equation}
\mathcal{O}, \mathcal{U}_1, \mathcal{U}_2,
\mathcal{O}(H), \mathcal{U}_1(H), \mathcal{U}_2(H), \dots,
\mathcal{O}((m-1)H), \mathcal{U}_1((m-1)H), \mathcal{U}_2((m-1)H),
\end{equation}
so that in particular \cref{question:first-semiorthogonal-embedding,question:r-semiorthogonal-embeddings,question:strongly-exceptional}
have a positive answer.
\end{proposition}
\begin{proof}
First, note that the Harder--Narasimhan types that appear are independent of $m$: the subdimension
vectors of $(2, 3)$ which admit $(3, -2)$-semistable representation continue
to do so if $m$ grows, so all $7$ possible nontrivial Harder--Narasimhan types occur.
They are listed in \cref{tab:m-kronecker-quiver-case}.

Now, since $\canonicalstability = m \cdot (3, -2)$, from
\cref{proposition:weights-universal-bundles,corollary:weight-uidualuj,proposition:weight-H}
we see that all the Teleman weights grow at most linearly in $m$.
On the other hand, from \cref{proposition:Teleman-bounds} we see that $\eta_{\tuple{d}^*}$
grows quadratically in $m$, so for $m$ large enough all the Teleman inequalities
\begin{align}
\max W (U^{\vee}_i \otimes U_{j} \otimes L(-\frac{s}{m}\canonicalstability), \tuple{d}^*) &< \eta_{\tuple{d}^*}, \label{equation:Teleman-inequalities-explicit-1}\\
\max W (U^{\vee}_i \otimes L(-\frac{s}{m}\canonicalstability), \tuple{d}^*) &< \eta_{\tuple{d}^*}, \label{equation:Teleman-inequalities-explicit-2}\\
\max W (U_{j} \otimes L(-\frac{s}{m}\canonicalstability), \tuple{d}^*) &< \eta_{\tuple{d}^*}, \text{ and} \label{equation:Teleman-inequalities-explicit-3}\\
\max W (L(-\frac{s}{m}\canonicalstability), \tuple{d}^*) &< \eta_{\tuple{d}^*} \label{equation:Teleman-inequalities-explicit-4}
\end{align}
will hold for the appropriate values of $s$. In fact, these are already satisfied when $m = 3$.
This suffices to show that all the higher cohomology groups in \cref{lemma:criterion-QB}
vanish.

By applying Serre duality as in \cref{lemma:serre-duality}, we can also conclude
that for all $i, j \in Q_0$ and for all $s \in 1, \dots, m - 1$, we have
\begin{align}
\HH^{0}(\modulispace, \mathcal{U}^{\vee}_{i} \otimes \mathcal{U}_{j} \otimes \mathcal{O}(-sH)) &= 0, \\
\HH^{0}(\modulispace, \mathcal{O}(-sH)) &= 0, \\
\HH^{0}(\modulispace, \mathcal{U}^{\vee}_{i} \otimes \mathcal{O}(-sH)) &= 0, \text{ and} \\
\HH^{0}(\modulispace, \mathcal{U}_{j} \otimes \mathcal{O}(-sH)) &= 0.
\end{align}
This shows that \eqref{equation:exceptional-collection-2-3-m-kronecker} is an exceptional collection.

To apply \cref{lemma:criterion:QC} now, we just notice that all the Teleman weights of $L(\frac{1}{m}\canonicalstability)$
are negative. This means that if the Teleman inequality is already satisfied on all the strata for
some vector bundle $F$, it will automatically be for $F \otimes L(\frac{1}{m}\canonicalstability)$.
By applying this to the terms in
the inequalities \eqref{equation:Teleman-inequalities-explicit-1}
to \eqref{equation:Teleman-inequalities-explicit-4},
we conclude.
\end{proof}

The first statement always yields a partial tilting object
\begin{equation}
\mathcal{T} \colonequals \bigoplus_{s = 0}^{m - 1} \mathcal{U} \otimes \mathcal{O}(sH).
\end{equation}

The second condition of the proposition, namely the fact that for all $i$ we have
\begin{equation}\label{equation:extra-vanishing}
\HH^{0}(\modulispace, \mathcal{U}^{\vee}_{i}) = 0,
\end{equation}
has been verified experimentally for values of $m$ up to $11$ in
\textsc{mkronecker-verification.jl} \cite{script},
and is expected to hold for all $m$.
When \eqref{equation:extra-vanishing} holds, the larger resulting partial tilting bundle is
\begin{equation}
\mathcal{T} \colonequals \bigoplus_{s = 0}^{m - 1}(\mathcal{O} \oplus \mathcal{U}) \otimes \mathcal{O}(sH).
\end{equation}

\begin{table}[H]
\centering
\begin{tabular}{ccccccc}
\toprule
$\tuple{d}^*$ &
$ \frac{1}{m} \cdot W(U^{\vee}_i \otimes U_j, \tuple{d}^*)$ &
$ \frac{1}{m} \cdot W(U_i(\tuple{a}), \tuple{d}^*)$ &
$ \frac{1}{m} \cdot W(L(\frac{\canonicalstability}{r}), \tuple{d}^*)$ &
$ \frac{1}{m} \cdot \eta_{\tuple{d}^*}$ \\
\midrule
((1, 1), (1, 2))        & $-\frac{5}{6}$ 0 $\frac{5}{6}$                                     & $-\frac{5}{6}$ 0                    & $-\frac{5}{6}$  & $\frac{5}{6}(2m - 3)$ \\
((2, 2), (0, 1))        & $-\frac{5}{2}$ 0 $\frac{5}{2}$                                     & $-5$ $-\frac{5}{2}$                   & $-5$              & $5(m - 1)$\\
((2, 1), (0, 2))        & $-\frac{10}{3}$ 0 $\frac{10}{3}$                                   & $-10$ $-\frac{20}{3}$                 & $-\frac{40}{3}$ & $\frac{10}{3}(4m - 2)$\\
((1, 0), (1, 3))        & $-\frac{15}{4}$ 0 $\frac{15}{4}$                                   & $-\frac{15}{2}$ $-\frac{15}{4}$     & $-\frac{45}{4}$ & $\frac{15}{4}(3m - 1)$\\
((1, 0), (1, 2), (0, 1))& $-5$ $-\frac{10}{3}$ $-\frac{5}{3}$ 0 $\frac{5}{3}$ $\frac{10}{3}$ 5 & $-10$ $-\frac{25}{3}$ -5              & $-\frac{40}{3}$ & $\frac{20}{3}(2m - 1)$\\
((1, 0), (1, 1), (0, 2))& $-5$ $-\frac{5}{2}$ 0 $\frac{5}{2}$ 5                                & $-\frac{25}{2}$ $-10$ $-\frac{15}{2}$ & $-\frac{35}{2}$ & $\frac{5}{2}(7m - 3)$\\
((2, 0), (0, 3))        & $-5$ 0 5                                                             & $-20$ $-15$                             & $-30$             & $30m$ \\
\bottomrule
\end{tabular}
\caption{Harder--Narasimhan types and non-normalized Teleman weights for \cref{proposition:m-kronecker-quivers}.}
\label{tab:m-kronecker-quiver-case}
\end{table}

Here we show directly by computation that
\cref{question:first-semiorthogonal-embedding,question:r-semiorthogonal-embeddings,question:strongly-exceptional}
have a positive answer for \cref{example:recurring-example}.

\begin{proposition}
The moduli space in \cref{example:recurring-example} admits a strongly exceptional collection
given by
\begin{equation}\label{equation:strongly-exceptional-sequence-recurring-example}
\mathcal{O}_{\modulispace}, \mathcal{U}_{1}, \mathcal{U}_{2},
\mathcal{O}_{\modulispace}(H), \mathcal{U}_{1}(H), \mathcal{U}_{2}(H),
\mathcal{O}_{\modulispace}(2H), \mathcal{U}_{1}(2H), \mathcal{U}_{2}(2H).
\end{equation}
In particular, \cref{question:first-semiorthogonal-embedding,question:r-semiorthogonal-embeddings,question:strongly-exceptional}
have a positive answer.
\end{proposition}
\begin{proof}
To show that the necessary vanishings outlined in
\eqref{equation:needed-UidualUjH}, \eqref{equation:needed-H},
\eqref{equation:needed-UiH} and \eqref{equation:needed-UidualH}
are satisfied, we run the script \textsc{verification.jl} \cite{script},
which verifies the necessary Teleman inequalities
and computes various Euler characteristics.
For the bundles $\mathcal{U}^{\vee}_i(-H)$ and $\mathcal{U}^{\vee}_i(-2H)$,
the Teleman inequality is not satisfied, so we still have to show that
\begin{align}
\HH^{\geq 1}(\modulispace, \mathcal{U}^{\vee}_i(- H)) &= 0, \text{ and}\\
\HH^{\geq 1}(\modulispace, \mathcal{U}^{\vee}_i(-2H)) &= 0.
\end{align}
From the output of \textsc{verification.jl} \cite{script} we know that for all $i \in Q_0$ and all $k \geq 0$,
\begin{align}
\HH^{k}(\modulispace, \mathcal{U}_i(-H)) &= 0, \text{ and}\\
\HH^{k}(\modulispace, \mathcal{U}_i(-2H)) &= 0,
\end{align}
so we can apply Serre duality as phrased in \cref{lemma:serre-duality} to conclude
that \eqref{equation:strongly-exceptional-sequence-recurring-example} is exceptional.

From \cref{tab:HN-recurring-example} we see that the weights of $H$ are negative on every
Harder--Narasimhan stratum.
Since the Teleman inequalities for
\eqref{equation:needed-UidualUjH}, \eqref{equation:needed-H},
\eqref{equation:needed-UiH} and \eqref{equation:needed-UidualH}
are all satisfied, and twisting by $H$ only lowers the Teleman weights,
we can conclude that for all $s \geq 1$,
\begin{align}
\HH^{\geq 1}(\modulispace, \mathcal{U}^{\vee}_{i} \otimes \mathcal{U}_{j} \otimes \mathcal{O}_{\modulispace}(sH)) &= 0, \\
\HH^{\geq 1}(\modulispace, \mathcal{O}_{\modulispace}(sH)) &= 0, \\
\HH^{\geq 1}(\modulispace, \mathcal{U}^{\vee}_{i} \otimes \mathcal{O}_{\modulispace}(sH)) &= 0, \text{ and}\\
\HH^{\geq 1}(\modulispace, \mathcal{U}_{j} \otimes \mathcal{O}_{\modulispace}(sH)) &= 0.
\end{align}
It remains to show that $\HH^{\geq 1}(\modulispace, \mathcal{U}_{j}) = 0$.
From \cref{tab:HN-recurring-example}, we see that on every Harder--Narasimhan stratum $\tuple{d}^*$
we have
\begin{equation}
\max(W(U_j(\tuple{a}), \tuple{d}^*)) < \eta_{\tuple{d}^*},
\end{equation}
so we conclude.
\end{proof}

\subsection{del Pezzo surfaces}\label{subsection:del-pezzos}
There exist 10 deformation families of del Pezzo surfaces. Of these,
4 have nontrivial infinitesimal deformations. Under \cref{assumptions}, the resulting
quiver moduli are infinitesimally rigid \cite[Corollary D]{rigidity-paper}, so up to
6 del Pezzo surfaces can be realised as quiver moduli, and it is known that,
in fact, all of them are. This is explained in \cite[p. 998]{MR4352662},
in \cite[Section 6]{king-note},
and stated below.

\begin{proposition}
Let $X$ be a rigid del Pezzo surface. There exist $Q, \tuple{d}$ such that
$X \simeq \modulispace^{\canonicalstability\stable}(Q, \tuple{d})$.
\end{proposition}
\begin{proof}
This statement can be proven case by case.
In each case shown in \cref{tab:del-pezzos},
the dimension of the Hochschild homology
is computed using \cite[Corollary 6.8]{MR1974891}.
Both the rank of each Picard group and the fact that these are Fano surfaces
are immediate consequences of \cite[Proposition 4.3]{MR4352662}.
Lastly, the index is in each case the $\gcd$ of $\canonicalstability$.
Comparing these invariants against the
well-known classification of del Pezzo surfaces (see e.g. \cite{fanography}),
we conclude that all $6$ of them are realised.
\end{proof}

\begin{figure}[H]
\label{del-Pezzo-quivers}
\centering
    \begin{tikzpicture}[baseline = -30pt, node distance=1.5cm]
    \node (1)                 {$\bullet$};
    \node (2) [below of = 1]  {$\bullet$};
    \node (l) [below =0.8cm of 2]  {$Q_1$};

    \draw[->] (1) edge (2);
    \draw[->, bend  left=15] (1) edge (2);
    \draw[->, bend right=15] (1) edge (2);
    \end{tikzpicture}\qquad~~
    \begin{tikzpicture}[baseline = -30pt, node distance=0.75cm]
    \node (1)                 {$\bullet$};
    \node (2) [below of = 1]  {$\bullet$};
    \node (3) [below of = 2]  {$\bullet$};
    \node (l) [below =0.8cm of 3]  {$Q_2$};

    \draw[->, bend  left=15] (1) edge (2);
    \draw[->, bend right=15] (1) edge (2);
    \draw[->, bend  left=15] (3) edge (2);
    \draw[->, bend right=15] (3) edge (2);
    \end{tikzpicture}\qquad~~
    \begin{tikzpicture}[baseline = -30pt, node distance=1.5cm]
    \node (1)                 {$\bullet$};
    \node (2) [right of = 1]  {$\bullet$};
    \node (3) [below of = 2]  {$\bullet$};
    \node (s) [below =0.8cm of 3]  {$Q_3$};

    \draw[->] (1) edge (2);
    \draw[->] (1) edge (3);
    \draw[->, bend  left=15] (2) edge (3);
    \draw[->, bend right=15] (2) edge (3);

    \end{tikzpicture}\qquad~~
    \begin{tikzpicture}[baseline = -30pt, node distance=1.5cm]
    \node (1)                 {$\bullet$};
    \node (3) [right of = 1]  {$\bullet$};
    \node (4) [below of = 1]  {$\bullet$};
    \node (2) [right of = 4]  {$\bullet$};
    \node (s) [below =0.8cm of 2] {$Q_4$};

    \draw[->] (1) edge (3);
    \draw[->] (1) edge (4);
    \draw[->] (2) edge (3);
    \draw[->] (2) edge (4);
    \draw[->] (3) edge (4);

    \end{tikzpicture}\qquad~~
    \begin{tikzpicture}[baseline = -30pt, node distance=0.75cm]
    \node (1)                 {$\bullet$};
    \node (2) [right of = 1]  {};
    \node (3) [right of = 2]  {$\bullet$};
    \node (4) [below of = 1]  {};
    \node (5) [right of = 4]  {$\bullet$};
    \node (6) [right of = 5]  {};
    \node (7) [below of = 4]  {$\bullet$};
    \node (8) [right of = 7]  {};
    \node (9) [right of = 8]  {$\bullet$};
    \node (l) [below =1.55cm of 5]  {$Q_5$};

    \draw[->] (1) edge (3);
    \draw[->] (1) edge (7);
    \draw[->] (5) edge (3);
    \draw[->] (5) edge (7);
    \draw[->] (9) edge (3);
    \draw[->] (9) edge (7);
    \end{tikzpicture}\qquad~~
    \begin{tikzpicture}[baseline = -30pt, node distance=0.5cm, column sep=1 .5cm]
    \node (1)                 {$\bullet$};
    \node (2) [below of = 1]  {$\bullet$};
    \node (3) [below of = 2]  {$\bullet$};
    \node (4) [below of = 3]  {$\bullet$};
    \node (5) [below of = 4]  {$\bullet$};
    \node (6) [right of = 3]  {};
    \node (7) [right of = 6]  {$\bullet$};
    \node (l) [below =1.4cm of 6]  {$Q_6$};

    \draw[->] (1) edge (7);
    \draw[->] (2) edge (7);
    \draw[->] (3) edge (7);
    \draw[->] (4) edge (7);
    \draw[->] (5) edge (7);
    \end{tikzpicture}
\caption{From left to right,
the quivers $Q_{1}$ to $Q_{6}$ realising the del Pezzo surfaces.
The corresponding surface to each quiver is described in \cref{tab:del-pezzos}.}
\end{figure}

We refer to the surfaces above as ``del Pezzo quiver moduli'', and we provide tilting
objects for all of them in the following proposition.

\begin{proposition}
For each del Pezzo quiver moduli $\modulispace$,
\cref{question:first-semiorthogonal-embedding} has a positive answer and yields
a full strongly exceptional collection of the form
\begin{equation}
\derived^{\bounded}(X) = \left\langle \mathcal{O}, \mathcal{U}_1, \dots, \mathcal{U}_{\#Q_0}\right\rangle.
\end{equation}
In particular, for each of them the vector bundle
\begin{equation}
\mathcal{T} \colonequals \mathcal{O}_{\modulispace} \oplus \mathcal{U}_{1} \dots \oplus \mathcal{U}_{\#Q_0}
\end{equation}
is tilting.
\end{proposition}
\begin{proof}
The proof of everything but fullness is a numerical verification of the fact that
all the necessary Euler characteristics vanish,
all the necessary Teleman inequalities are satisfied, and that whenever the Teleman
inequality is not satisfied directly, Serre duality can be applied as in \cref{lemma:serre-duality}.
All the computations and Serre duality checks are performed in \textsc{verification.jl} \cite{script}.

The fact that exceptional collections of maximal length are full is a consequence of
\cite[Theorem 6.35]{MR4581971}.
\end{proof}

\begin{remark}\label{example:p2-p1-times-p1}
For $\PP^2$ and $\PP^1 \times \PP^1$, \cref{question:r-semiorthogonal-embeddings} has
a negative answer. This is an example of the Hochschild
homology behaviour discussed in \cref{remark:hochschild-homology}:
for $\PP^2$, respectively $\PP^1 \times \PP^1$,
\cref{question:r-semiorthogonal-embeddings}
predicts an exceptional collection of length 9, respectively 8,
but the maximum length is 3, respectively 4.
\end{remark}

\begin{table}[H]
\centering
\begin{tabular}{cccccc}
\toprule
del Pezzo &
$\modulispace^{\theta_{\mathrm{can}}}(Q, \tuple{d})$ &
$\dim \Hochschild_{0}(\modulispace)$ &
$\rk(\Pic(\modulispace))$ &
$r$ &
Full exc. collection. \\
\midrule
$\PP^2$                 & $\modulispace^{\theta_{\mathrm{can}}}(Q_{1}, \tuple{1})$          & $3$ & $1$ & $3$ & $\mathcal{O}, \mathcal{U}_{1}, \mathcal{U}_{2}$ \\
$\PP^1 \times \PP^1$    & $\modulispace^{\theta_{\mathrm{can}}}(Q_{2}, \tuple{1})$          & $4$ & $2$ & $2$ & $\mathcal{O}, \mathcal{U}_{1}, \mathcal{U}_{2}, \mathcal{U}_{3}$ \\
$\mathrm{Bl}_{1}(\PP^2)$& $\modulispace^{\theta_{\mathrm{can}}}(Q_{3}, \tuple{1})$          & $4$ & $2$ & $1$ & $\mathcal{O}, \mathcal{U}_{1}, \mathcal{U}_{2}, \mathcal{U}_{3}$ \\
$\mathrm{Bl}_{2}(\PP^2)$& $\modulispace^{\theta_{\mathrm{can}}}(Q_{4}, \tuple{1})$          & $5$ & $3$ & $1$ & $\mathcal{O}, \mathcal{U}_{1}, \mathcal{U}_{2}, \mathcal{U}_{3}, \mathcal{U}_{4} $ \\
$\mathrm{Bl}_{3}(\PP^2)$& $\modulispace^{\theta_{\mathrm{can}}}(Q_{5}, \tuple{1})$          & $6$ & $4$ & $1$ & $\mathcal{O}, \mathcal{U}_{1}, \mathcal{U}_{2}, \mathcal{U}_{3}, \mathcal{U}_{4},  \mathcal{U}_{5} $ \\
$\mathrm{Bl}_{4}(\PP^2)$& $\modulispace^{\theta_{\mathrm{can}}}(Q_{6}, (1, \dots, 1, 2))$   & $7$ & $5$ & $1$ & $\mathcal{O}, \mathcal{U}_{1}, \mathcal{U}_{2}, \mathcal{U}_{3}, \mathcal{U}_{4},  \mathcal{U}_{5},  \mathcal{U}_{6} $ \\
\bottomrule
\end{tabular}
\caption{Rigid del Pezzo surfaces as quiver moduli and their exceptional collections.}
\label{tab:del-pezzos}
\end{table}

We now present a higher-dimensional example, namely a Fano $5$-fold.
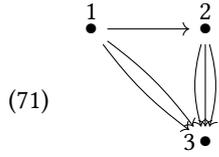
\begin{example}
\begin{equation}\label{fig:5fold-3-vertex-quiver}
\begin{tikzpicture}[baseline = -30, node distance = 1.5cm]
\node (1)                 {$\bullet$};
\node (2) [right of = 1]  {$\bullet$};
\node (3) [below of = 2]  {$\bullet$};

\draw (1) node[above] {1};
\draw (2) node[above] {2};
\draw (3) node[left]  {3};

\draw[->]                  (1) edge (2);
\draw[->, bend  left = 8] (1) edge (3);
\draw[->, bend right = 8] (1) edge (3);

\draw[->, bend  left = 15]  (2) edge (3);
\draw[->, bend  left = 8]   (2) edge (3);
\draw[->, bend right = 8]   (2) edge (3);
\draw[->, bend right = 15]  (2) edge (3);

\end{tikzpicture}
\end{equation}
Let $Q$ be the $3$-vertex quiver in \eqref{fig:5fold-3-vertex-quiver}.
Let $\tuple{d} = (1, 1, 1)$, and let $\tuple{a} = (0, 2, -1)$.
The resulting moduli space $\modulispace$ is a Fano $5$-fold with index $3$,
Picard rank $2$ and $\dim \Hochschild_{0}(\modulispace) = 12$.
We can compute the degree of $-K_{\modulispace}$ in the Chow ring, and it turns out that
$(-K_{\modulispace})^5 = 6318$.

Running \textsc{verification.jl} \cite{script}, we see that \cref{question:r-semiorthogonal-embeddings,question:strongly-exceptional}
not only give a strongly exceptional collection, but said collection is also of maximal length.
\end{example}

\subsection{Non-examples}
As mentioned in \cref{section:introduction},
cases are known where
\cref{question:first-semiorthogonal-embedding,question:r-semiorthogonal-embeddings}
have negative answers.
The only occurring quiver moduli of dimension $1$, discussed below, is such an example,
while all the occurring surfaces have been discussed in \cref{subsection:del-pezzos}.

\begin{example}\label{example:p1}
\begin{equation}\label{equation:Kronecker-quiver}
\begin{tikzpicture}[baseline = -3pt, node distance = 1.5cm]
\node (1)                {$\bullet$};
\node (2) [right of = 1] {$\bullet$};

\draw (1) node[above] {1};
\draw (2) node[above] {2};

\draw[->, bend left = 15] (1) edge (2);
\draw[->, bend right= 15] (1) edge (2);
\end{tikzpicture}
\end{equation}
Let $Q$ be the Kronecker quiver as in \eqref{equation:Kronecker-quiver},
with two vertices $1$ and $2$ and two arrows from $1$ to $2$.
Let $\tuple{d} = (1, 1)$.
It is easy to show directly that \cref{assumptions} is satisfied,
and the resulting quiver moduli $\modulispace^{\canonicalstability\stable}(Q, \tuple{d})$
is isomorphic to $\PP^1$.
We know that the index of $\PP^1$ is 2, and that $\dim \Hochschild_{0}(\PP^1) = 2$.
Therefore, \cref{question:first-semiorthogonal-embedding}, which predicts an exceptional
collection of $3$ objects, must have a negative answer.
This is the smallest possible example of the pathological behaviour described in
\cref{remark:hochschild-homology}.
\end{example}

We exhibit now a case in dimension $3$ where
\cref{question:r-semiorthogonal-embeddings} has a negative answer.

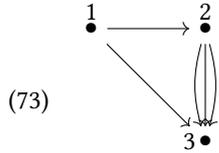
\begin{example}
\begin{equation}\label{fig:3-vertex-quiver}
\begin{tikzpicture}[baseline = -30, node distance = 1.5cm]
\node (1)                 {$\bullet$};
\node (2) [right of = 1]  {$\bullet$};
\node (3) [below of = 2]  {$\bullet$};

\draw (1) node[above] {1};
\draw (2) node[above] {2};
\draw (3) node[left]  {3};

\draw[->] (1) edge (2);
\draw[->] (1) edge (3);

\draw[->, bend  left = 15]  (2) edge (3);
\draw[->]                   (2) edge (3);
\draw[->, bend right = 15]  (2) edge (3);

\end{tikzpicture}
\end{equation}
Let now $Q$ be the three-vertex quiver in \eqref{fig:3-vertex-quiver}.
Let $\tuple{d} = (1, 1, 1)$ and $\tuple{a} = (1, 1, -1)$.
The resulting moduli space $\modulispace \colonequals\modulispace^{\canonicalstability\stable}(Q, \tuple{d})$
is a Fano threefold of index $2$, with $\dim \Hochschild_{0}(\modulispace) = 6$ and
with $(-K_{\modulispace})^3 = 56$.
These invariants suffice to conclude that $\modulispace$ is the Fano variety 2-35
\cite[\href{https://www.fanography.info/2-35}{2-35}]{fanography}.

Running \textsc{verification.jl} \cite{script}, we see that
\cref{question:first-semiorthogonal-embedding,question:strongly-exceptional}
hold, so that the sequence we obtain is strongly exceptional. \cref{question:r-semiorthogonal-embeddings}
however can't hold, for the reason described in \cref{remark:hochschild-homology}.
\end{example}

\subsection{The proof of \texorpdfstring{\cref{result:weaker-theorem}}{Theorem D}}\label{subsection:a-weaker-result}
As stated in the introduction, it is possible to give a sufficient criterion
for a weaker version
of \cref{question:r-semiorthogonal-embeddings,question:strongly-exceptional} that is
independent of the choice of linearisation of \cref{remark:choice-linearisation} to hold.
In this section we collect the necessary results to prove~\cref{result:weaker-theorem}.

For each Harder--Narasimhan type $\tuple{d}^* = (\tuple{d}^1,\dots,\tuple{d}^{\ell})$
and for each $s = 1,\dots,\ell$,
let~$k_{\tuple{d}^*, s} \colonequals c \cdot \mu_{\canonicalstability}(\tuple{d}^s)$,
where~$c$ is the smallest positive integer
for which~$k_{\tuple{d}^*, s}$ is integer for every~$s$.
Let~$t_{\tuple{d}^*}$ be the largest integer ${t}$ for which we have the inequality
\begin{equation}\label{equation:required-inequality-weaker-theorem}
k_{\tuple{d}^*, 1} - k_{\tuple{d}^*, {\ell}} +
\frac{{t}}{r}\canonicalstability\cdot\left(\sum_{s = 1}^{\ell}k_{\tuple{d}^*, s}\tuple{d}^s\right) <
\eta_{\tuple{d}^*}.
\end{equation}

From \cite[Proposition~4.2]{rigidity-paper} we know that $\min_{\tuple{d}^*}\{t_{\tuple{d}^*}\}$
must be nonnegative,
and from~\cite[Theorem~A]{vector-fields-paper} we know that the same quantity
is bounded above by~$r - 1$.
We will now prove that, when this value is maximal, \cref{result:weaker-theorem} holds.

\begin{assumption}\label{assumption:t-large}
We have $\min_{\tuple{d}^*}\{t_{\tuple{d}^*}\} = r - 1$.
\end{assumption}

\begin{theorem}\label{theorem:weaker-sod}
Let $Q, \tuple{d}$ and the corresponding $\canonicalstability$
satisfy \cref{assumptions}.
If \cref{assumption:t-large} is satisfied,
there is a semiorthogonal decomposition of $\derived^{\bounded}(\modulispace)$ given by
\begin{equation}\label{equation:smaller-sod}
\derived^{\bounded}(\modulispace) =
\left\langle \Phi_{\mathcal{U}}(\derived^{\bounded}(Q)), \dots, \Phi_{\mathcal{U}((r-1)H)}(\derived^{\bounded}(Q));\dots \right\rangle.
\end{equation}
\end{theorem}
\begin{proof}

To verify that the terms in \eqref{equation:smaller-sod} are left orthogonal in the given order,
we only need to verify the vanishings in~\eqref{equation:needed-UidualUjH}.
As all the Teleman weights of ${L}(\frac{1}{r}\canonicalstability)$ are negative,
in order to apply \cref{corollary:teleman-vanishing}
to~$U^{\vee}_i \otimes U_j \otimes {L}(-\frac{{t}}{r}\canonicalstability)$
for every~${t} = 1,\dots,r-1$ it is enough to verify, for every Harder--Narasimhan type, that
\begin{equation}
k_{\tuple{d}^*, 1} - k_{\tuple{d}^*, {\ell}} + \frac{r - 1}{r}\canonicalstability\cdot\left(\sum_{s = 1}^{\ell}k_{\tuple{d}^*, s}\tuple{d}^s\right) < \eta_{\tuple{d}^*}.
\end{equation}
\cref{assumption:t-large} and Serre duality as stated in~\cref{lemma:serre-duality}
allow to conclude that for all $k \in \mathbb{Z}$,
for all $i, j \in Q_{0}$ and all ${t} = 1,\dots,r - 1$,
we have~$\HH^{k}(\modulispace, \mathcal{U}^{\vee}_{i}\otimes\mathcal{U}_{j}\otimes\mathcal{O}(-{t}H)) = 0$.
\end{proof}

\begin{proposition}\label{proposition:weaker-strong-exceptional}
In the same setting of \cref{theorem:weaker-sod},
the quiver moduli $\modulispace$ admits a strong exceptional collection given by
\begin{equation}\label{equation:weaker-theorem-strong-exceptional-collection}
\mathcal{U}_{1},\dots,\mathcal{U}_{n},\dots,\mathcal{U}_{1}((r-1)H),\dots,\mathcal{U}_{n}((r-1)H).
\end{equation}
\end{proposition}
\begin{proof}
We must show that for all $t = 0,\dots,r-1$,
we have~$\HH^{\geq 1}(\modulispace,\mathcal{U}^{\vee}_{i}\otimes\mathcal{U}_{j}\otimes\mathcal{O}(tH)) = 0$.
The case $t = 0$ is given, under \cref{assumptions}, by
applying \cref{corollary:teleman-vanishing} using~\cite[Proposition~4.2]{rigidity-paper}.
For $t > 0$, \cref{corollary:teleman-vanishing} can be applied verbatim, as
the Teleman weights of ${L}(\frac{1}{r}\canonicalstability)$ are always negative.
\end{proof}
To conclude, \cref{theorem:weaker-sod,proposition:weaker-strong-exceptional} together
prove \cref{result:weaker-theorem}.
\printbibliography

\emph{Gianni Petrella} \url{gianni.petrella@uni.lu} \\
Universit\'e de Luxembourg, 6, Avenue de la Fonte, L-4364 Esch-sur-Alzette, Luxembourg

\end{document}